\documentclass[10pt]{article}%
\usepackage{amsfonts}
\usepackage{amsmath}
\usepackage{amssymb}
\usepackage{graphicx}
\usepackage[left=2.5cm,top=3.0cm,right=2.5cm,bottom=2.5 cm]{geometry}%
\usepackage{hyperref}
\usepackage{xcolor}
\providecommand{\U}[1]{\protect\rule{.1in}{.1in}}
\newtheorem{theorem}{Theorem}

\newtheorem{proposition}[theorem]{Proposition}
\newtheorem{remark}[theorem]{Remark}

\newenvironment{proof}[1][Proof]{\noindent\textbf{#1.} }{\ \rule{0.5em}{0.5em}}
\begin{document}

\title{\textbf{Blow-up of solutions to the Keller-Segel model with tensorial flux in
high dimensions}}
\author{Valeria Cuentas, Elio Espejo \thanks{Corresponding author. \newline E-mail
addresses: valeria.cuentasrodriguez@nottingham.edu.cn (V. Cuentas),
elio-eduardo-espejo.arenas@nottingham.edu.cn(E. Espejo),
suzuki@sigmath.es.osaka-u.ac.jp (T. Suzuki)}, \ and Takashi Suzuki}
\date{}
\maketitle

\begin{abstract}
Over the course of the last decade, there has been a significant level of
interest in the analysis of Keller-Segel models incorporating tensorial flux.
Despite this interest, the question of whether finite-time blowup solutions
exist remains a topic of ongoing research. Our study provides evidence that
solutions of this nature are indeed possible in dimensions $n\geq3,$ when
utilizing a tensorial flux expressed in the form of $A\nabla v$, where $A$
denotes a matrix with constant components.\bigskip


\end{abstract}

\noindent
\textcolor{red}{
For access to the published version of this paper in \textit{Applied Mathematics Letters}, Volume 154, August 2024, please refer to the following link: 
\href{https://doi.org/10.1016/j.aml.2024.109090}{Journal version available here}.
}

\vspace{0.5cm}

\section{Introduction}

Chemotaxis is an intriguing biological phenomenon that plays a crucial role in
enabling the aggregation and distribution of various species. It is a process
that involves the movement of cells or organisms towards a chemical gradient,
which is a concentration of molecules that stimulates the cells or organisms
to move in a particular direction. Chemotaxis is an essential mechanism in
many biological processes, including the immune response, wound healing, and
embryonic development. It is also a critical factor in the behavior of
microorganisms, such as bacteria, which use chemotaxis to locate nutrients and
avoid toxins. Thus, the study of chemotaxis is essential to understanding the
behavior and interactions of living organisms at the molecular level. This
process involves the movement of organisms in response to a concentration
gradient of chemicals. The model developed by Keller and Segel is widely
recognized as a seminal contribution to the field of chemotaxis. It provides a
mathematical framework for understanding the mechanisms underlying this
complex biological process(e.g \cite{KellerSegel}). This model can be
simplified by%
\begin{equation}
u_{t}=\Delta u-\chi\nabla\cdot(u\nabla v),\text{ \ \ \ and \ \ }\varepsilon
v_{t}=\Delta v-v+u, \label{KS1}%
\end{equation}
\bigskip where $u(x,t)$ denotes the density and $v(x,t)$ the chemical
concentration at a given point $x$ and time $t.$

The model (\ref{KS1}) can exhibit interesting variations, particularly when
the migration is not parallel to the signal gradient. A notable example of
this phenomenon is exhibited by peritrichously flagellated bacteria when
swimming in close proximity to surfaces. In such cases, the density of
bacteria evolves according to the equation
\[
u_{t}=\Delta u-\nabla\cdot(uA(x,u,v)\nabla v),
\]
where $x\in%
\mathbb{R}
^{n},$ $t>0$ and $A(x,u,v)$ is a $n\times n$ matrix. Over the last decade,
several studies have been conducted on the global existence and asymptotic
behavior of this type of model with tensorial chemotaxis, including references
to (\cite{EspejoWu,WinklerMatrix,Cuentas}). Despite this progress, it remains
unclear whether solutions may experience blow up in finite time when the
chemoattractant is produced by the cells.\bigskip

We aim to prove the possibility of having solutions blowing-up in a finite
time for system%
\begin{equation}%
\begin{array}
[c]{cc}%
\partial_{t}u=\Delta u-\chi\nabla\cdot(uA\nabla v), & x\in\mathbb{R}%
^{n},t>0,\\
-\Delta v=u,\text{ }v(x,t)=\frac{1}{n(n-2)\left\vert B_{1}(0)\right\vert }\int
u(y,t)\left\vert x-y\right\vert ^{2-n}dy & x\in\mathbb{R}^{n},t>0,\\
u(x,0)=u_{0}(x)\geq0, & x\in\mathbb{R}^{n},
\end{array}
\label{Problem Matrix A nxn}%
\end{equation}
where $A:=(a_{ij})_{i,j=1,n}\in M_{n}(%
\mathbb{R}
)$ represents a nonsingular $n\times n$ matrix with constant components
satisfying $x^{T}\left(  \left(  AA^{T}\right)  ^{1/2}\right)  ^{-1}Ax>0$ for
all non-zero $x\in\mathbb{R}^{n}.$ Here the symbol $\sqrt{AA^{T}}$ stands for
the positive-definite square root of the matrix $AA^{T},$ whose existence and
uniqueness is well-established in mathematics (c.f. \cite[Corollary
7.3.3]{Roger}). Examples of matrices satisfying this hypothesis include the
set of positive-definite matrices and, in the three in the three-dimensional
case, orthogonal matrices of the form%
\[
A=\left(
\begin{array}
[c]{ccc}%
\cos\alpha & -\sin\alpha & 0\\
\sin\alpha & \cos\alpha & 0\\
0 & 0 & 1
\end{array}
\right)  ,
\]
where $\alpha\in(-\pi/2,\pi/2)$. Our approach to proving blow-up involves
decomposing matrix $A$ into its polar components and employing a modified
version of the second moments technique. In contrast to the nontensorial
Keller-Segel model, where the evolution of $\int_{%
\mathbb{R}
^{n}}u(x,t)\left\vert x\right\vert ^{2}dx$ is fundamental, we reveal that the
tensorial attraction makes $\int_{%
\mathbb{R}
^{n}}u(x,t)(x^{T}Bx)dx$ crucial, where the matrix $B$, with constant
component, is meticulously chosen to yield the desired outcome of blow-up.

\section{\textbf{Local existence, regularity, uniqueness, mass conservation
and non-negativity for arbitrary matrices}}

\begin{proposition}
\label{LocalExistence nxn}Let $n\geq3$ and $A\in M_{n}(%
\mathbb{R}
),$ and suppose that the initial data $u_{0}\in BUC(\mathbb{R}^{n})\cap
L^{1}(\mathbb{R}^{n})$ is non-negative. Then, there exist $T_{\max}\in\left(
0,+\infty\right]  $ and a non-negative%
\[
u\in C^{0}\left(  \left[  0,T_{\max}\right)  ;BUC(\mathbb{R}^{n})\right)  \cap
C^{0}\left(  \left[  0,T_{\max}\right)  ;L^{1}(\mathbb{R}^{n})\right)  \cap
C^{\infty}\left(  \mathbb{R}^{n}\times\left(  0,T_{\max}\right)  \right)  ,
\]
such that writing $v(\cdot,t)=\mathbf{K_n}(x)\ast u(\cdot,t),t\in\left(
0,T_{\max}\right)  ,$ with $\mathbf{K_n}(x):=\frac{1}{n(n-2)\left\vert
B_{1}(0)\right\vert }$ $|x|^{2-n},x\in\mathbb{R}^{n}\backslash\{0\}.$ we
obtain $v\in C^{\infty}\left(  \mathbb{R}^{n}\times\left(  0,T_{\max}\right)
\right)  ,$ $\nabla v\in L_{loc}^{\infty}(\left[  0,T_{\max}\right)
;L^{\infty}$ $(\mathbb{R}^{n};\mathbb{R}^{n})),$ and that $\left(  u,v\right)
$ forms a classical solution of (\ref{Problem Matrix A nxn}) in $\mathbb{R}%
^{n}\times\left(  0,T_{\max}\right).$ We also have the next extensibility
criterion,%
\[%
\begin{array}
[c]{c}%
\text{if }T_{\max}<+\infty,\text{ then both }\limsup_{t\rightarrow T_{\max}%
}\left\Vert u(\cdot,t)\right\Vert _{L^{\infty}(\mathbb{R}^{n})}=+\infty\\
\text{and }\limsup_{t\rightarrow T_{\max}}\left\Vert \nabla v(\cdot
,t))\right\Vert _{L^{\infty}(\mathbb{R}^{n})}=+\infty.
\end{array}
\]
This solution is uniquely determined in the sense that if $T\in\left(
0,T_{\max}\right)  ,$ and if $\left(  \widehat{u},\widehat{v}\right)  $ is a
classical solution of (\ref{Problem Matrix A nxn}) in $\mathbb{R}^{n}%
\times\left(  0,T_{\max}\right)  $ fulfilling $\widehat{u}\in C^{0}(\left[
0,T\right]  ;BUC$ $(\mathbb{R}^{n}))\cap C^{0}\left(  \left[  0,T\right]
;L^{1}(\mathbb{R}^{n})\right)  \cap C^{2,1}\left(  \mathbb{R}^{n}\times\left(
0,T\right)  \right)  $ and $\widehat{v}\in C^{2,0}\left(  \mathbb{R}^{n}%
\times\left(  0,T\right)  \right)  $ as well as $\nabla\widehat{v}\in
L^{\infty}\left(  \mathbb{R}^{n}\times\left(  0,T\right)  ;\mathbb{R}%
^{n}\right)  ,$ then $\widehat{u}\equiv u$ in $\mathbb{R}^{n}\times\left(
0,T\right).$ Moreover,
\begin{equation}
\int_{%
\mathbb{R}
^{n}}u(x,t)dx=\int_{%
\mathbb{R}
^{n}}u_{0}dx=:M\text{ for all }t\in\left(  0,T_{\max}\right).
\label{MassConservationMatrixAnxn}%
\end{equation}

\end{proposition}

\begin{proof}
See \cite[Proposition 1.1.]{WinklerMatrix}.
\end{proof}

\section{\textbf{Blow-up}}

Our methodology to establish blow-up in high dimensions hinges upon the
technique recently proposed for the analysis of blow-up for two-dimensional
Keller-Segel type systems with tensorial flux, cf. \cite{Cuentas}. This
methodology can be outlined in two key steps: firstly, leveraging the polar
decomposition of the tensor $A$ and secondly, examining the evolution of the
quantity $\int u(x^{T}Bx)dx$ using a strategically chosen matrix $B$ with
constant components.

\begin{theorem}
[Blow-up]\label{blow-up nxn}Given $n\geq3,$ consider a non-negative classical
solution $u$ of system (\ref{Problem Matrix A nxn}) with non-negative initial
data $u_{0}\in BUC(\mathbb{R}^{n})\cap L^{1}(\mathbb{R}^{n})$ and
$u_{0}\left\vert x\right\vert ^{2}\in L^{1}(\mathbb{R}^{n}).$ Suppose also
that $A\in M_{n}(%
\mathbb{R}
)$ is a nonsingular matrix with constant components satisfying%
\begin{equation}
x^{T}\left(  \left(  AA^{T}\right)  ^{1/2}\right)  ^{-1}Ax>0\text{ for all
non-zero }x\in\mathbb{R}^{n}. \label{Hypothesis}%
\end{equation}
Let $[0,T_{\max})$ be the maximal interval of local existence of the solution
guaranteed by Proposition \ref{LocalExistence nxn}. If the integral $m_0:=\int_{%
\mathbb{R}
^{n}}u_{0}\left\vert x\right\vert ^{2}dx$ is small enough compared to the mass $M$, more precisely, if for a constant $C_{Bl}:=C(A,\chi,n)>0$
\begin{equation}
    \int_{\mathbb{R}^{n}}u_{0}\left\vert x\right\vert ^{2}dx
    \leq C_{Bl}M^{\frac{n}{n-2}}, \label{SmallMoment}
\end{equation} then $T_{\max}<+\infty.$
\end{theorem}

\begin{proof}
To facilitate the presentation, we conduct a formal calculation of the
evolution of moments, assuming the solution $u$ is suitably regular and decay
sufficiently fast at infinity. We start by decomposing the nonsingular matrix
$A$ into the polar form $A=PU,$ where $P=(p_{ij})_{i,j=1,n}:=\left(
AA^{T}\right)  ^{1/2}$ is positive-definite and $U:=P^{-1}A$ is orthogonal
(cf. \cite[Corollary 7.3.3.]{Roger}). Next, we proceed to modify the
second-moment blow-up technique by multiplying the equation for the cell
density $u$ by the quadratic form $x\cdot Bx,$ where $B$ is a positive
definite matrix to be determined. Integrating the product, we obtain%
\[
\frac{d}{dt}\int_{%
\mathbb{R}
^{n}}u\left(  x\cdot Bx\right)  dx=\int_{%
\mathbb{R}
^{n}}\left(  x\cdot Bx\right)  \Delta udx-\chi\int_{%
\mathbb{R}
^{n}}\left(  x\cdot Bx\right)  \nabla\cdot\left(  uPU\nabla v\right)  dx.
\]
Integration by parts leads to%
\[
\frac{d}{dt}\int_{%
\mathbb{R}
^{n}}u\left(  x\cdot Bx\right)  dx=\int_{%
\mathbb{R}
^{n}}\Delta\left(  x\cdot Bx\right)  udx+\chi\int_{%
\mathbb{R}
^{n}}\nabla\left(  x\cdot Bx\right)  \left(  uPU\nabla v\right)  dx.
\]
Considering the symmetry of the matrix $B,$ the formula $\nabla\left(  x\cdot
Bx\right)  =2Bx$ holds$,$ and therefore%
\[
\frac{d}{dt}\int_{%
\mathbb{R}
^{n}}u\left(  x\cdot Bx\right)  dx=\int_{%
\mathbb{R}
^{n}}\Delta\left(  x\cdot Bx\right)  udx+\chi\int_{%
\mathbb{R}
^{n}}2Bx\cdot\left(  uPU\nabla v\right)  dx.
\]
Utilizing again the symmetry of the matrix $B$, the last integral can be
rewritten as%
\[
\int_{%
\mathbb{R}
^{n}}2Bx\cdot\left(  uPU\nabla v\right)  dx=2\int_{%
\mathbb{R}
^{n}}x\cdot\left(  BPU\nabla v\right)  udx.
\]
Consequently, we choose $B=P^{-1}$ to simplify the subsequent calculations.
This leads to
\[
\frac{d}{dt}\int_{%
\mathbb{R}
^{n}}u\left(  x\cdot P^{-1}x\right)  dx=\int_{%
\mathbb{R}
^{n}}\Delta\left(  x\cdot P^{-1}x\right)  udx+2\chi\int_{%
\mathbb{R}
^{n}}x\cdot\left(  U\nabla v\right)  udx.
\]
Direct computations yield $\Delta\left(  x\cdot P^{-1}x\right)  =2Tr(P^{-1}).$
Thus%
\[
\frac{d}{dt}\int_{%
\mathbb{R}
^{n}}u\left(  x\cdot P^{-1}x\right)  dx=2Tr(P^{-1})\int_{%
\mathbb{R}
^{n}}udx+2\chi\int_{%
\mathbb{R}
^{n}}x\cdot\left(  U\nabla(\mathbf{K_n}\ast u)\right)  udx.
\]
This expression can be further simplified using the mass conservation property
(\ref{MassConservationMatrixAnxn}) to obtain%
\begin{align*}
\frac{d}{dt}\int_{%
\mathbb{R}
^{n}}u\left(  x\cdot P^{-1}x\right)  dx  &  =2Tr(P^{-1})\int_{%
\mathbb{R}
^{n}}u_{0}dx+2\chi\int_{%
\mathbb{R}
^{n}}x\cdot\left(  U\nabla(\mathbf{K_n}\ast u)\right)  udx\\
&  =2Tr(P^{-1})M+2\chi\int_{%
\mathbb{R}
^{n}}x\cdot\left(  U\nabla(\mathbf{K_n}\ast u)\right)  udx.
\end{align*}
We now proceed to show that the orthogonality of matrix $U$ allows for a
significant reduction of the integral $\int_{%
\mathbb{R}
^{n}}x\cdot\left(  uU\nabla(\mathbf{K_n}\ast u)\right)  dx.$ First, we
explicitly write the convolution $\nabla\left(  \mathbf{K_n}\ast u\right)  $ to
get%
\begin{align}
&  \frac{d}{dt}\int_{%
\mathbb{R}
^{n}}u\left(  x\cdot P^{-1}x\right)  dx\nonumber\\
&  =2Tr(P^{-1})M+2\chi\int_{%
\mathbb{R}
^{n}}x\cdot\left(  U\nabla(\mathbf{K_n}\ast u)\right)  udx\nonumber\\
&  =2Tr(P^{-1})M+2\chi\int_{%
\mathbb{R}
^{n}}x\cdot U\left(  \frac{-1}{n\left\vert B_{1}(0)\right\vert }\int_{%
\mathbb{R}
^{n}}\frac{x-y}{\left\vert x-y\right\vert ^{n}}u(y,t)dy\right)
u(x,t)dx\nonumber\\
&  =2Tr(P^{-1})M-\frac{2\chi}{n\left\vert B_{1}(0)\right\vert }\int_{%
\mathbb{R}
^{n}}\int_{%
\mathbb{R}
^{n}}\left(  x\cdot U\frac{x-y}{\left\vert x-y\right\vert ^{n}}%
u(x,t)u(y,t)dy\right)  dxdy. \label{lolan}%
\end{align}
We interchange $x$ and $y$ in the last integral to obtain
\begin{align*}
&  \int_{%
\mathbb{R}
^{n}}\int_{%
\mathbb{R}
^{n}}\left(  x\cdot U\frac{x-y}{\left\vert x-y\right\vert ^{n}}%
u(x,t)u(y,t)dy\right)  dxdy\\
&  =-\int_{%
\mathbb{R}
^{n}}\int_{%
\mathbb{R}
^{n}}\left(  y\cdot U\frac{x-y}{\left\vert x-y\right\vert ^{n}}%
u(x,t)u(y,t)dy\right)  dx,
\end{align*}
which implies%
\begin{align*}
&  \int_{%
\mathbb{R}
^{n}}\int_{%
\mathbb{R}
^{n}}\left(  x\cdot U\frac{x-y}{\left\vert x-y\right\vert ^{n}}%
u(x,t)u(y,t)dy\right)  dxdy\\
&  =\frac{1}{2}\int_{%
\mathbb{R}
^{n}\times%
\mathbb{R}
^{n}}\left(  \left(  x-y\right)  \cdot U\frac{x-y}{\left\vert x-y\right\vert
^{n}}u(x,t)u(y,t)dy\right)  dx.
\end{align*}
Thus, the identity (\ref{lolan}) reduces to%
\begin{align*}
&  \frac{d}{dt}\int_{%
\mathbb{R}
^{n}}u\left(  x\cdot P^{-1}x\right)  dx\\
&  =2Tr(P^{-1})M-\frac{\chi}{n\left\vert B_{1}(0)\right\vert }\int_{%
\mathbb{R}
^{n}\times%
\mathbb{R}
^{n}}\left(  \left(  x-y\right)  \cdot U\frac{x-y}{\left\vert x-y\right\vert
^{n}}u(x,t)u(y,t)dy\right)  dx.
\end{align*}
Next, we observe that since $U$ is an orthogonal matrix, there is an
orthogonal matrix $Q$ and a block diagonal matrix $D$ such that%
\begin{equation}
QUQ^{T}=D=%
\begin{pmatrix}
R_{1} &  &  &  &  & \\
& \ddots &  &  & \mathbf{0} & \\
&  & R_{k} &  &  & \\
&  &  & \lambda_{1} &  & \\
& \mathbf{0} &  &  & \ddots & \\
&  &  &  &  & \lambda_{p}%
\end{pmatrix}
, \label{CF}%
\end{equation}
where all the $R_{j}$ represent a $2\times2$ rotation matrix (cf.
\cite[Corollary 2.5.14. (c)]{Roger})$,$ that is a matrix of the form%
\[
R_{j}=\left(
\begin{array}
[c]{cc}%
\cos\alpha_{j} & -\sin\alpha_{j}\\
\sin\alpha_{j} & \cos\alpha_{j}%
\end{array}
\right)  ,\text{ where }\alpha_{j}\in(-\pi,\pi],
\]
and each $\lambda_{j}$ can take solely the values $1$ or $-1.$ Moreover, the
hypothesis that $0<x^{T}\left(  \left(  AA^{T}\right)  ^{1/2}\right)
^{-1}Ax=x^{T}P^{-1}Ax=x^{T}Ux$ for all non-zero $x\in%
\mathbb{R}
^{n}$, readily implies that $\lambda_{i}=1,i=1,\ldots,p,$ and $\cos\alpha
_{j}>0,j=1,\ldots,k.$Therefore, for any $x\in%
\mathbb{R}
^{n},$%
\begin{align*}
x^{T}Ux  &  =x^{T}Q^{T}DQx=(Qx)^{T}D(Qx)=(Qx)^{T}\left(  \frac{1}{2}\left(
D+D^{T}\right)  \right)  (Qx)\\
&  =(Qx)^{T}%
\begin{pmatrix}
\cos\alpha_{1} &  &  &  &  &  &  & \\
& \cos\alpha_{1} &  &  &  &  & \mathbf{0} & \\
&  & \ddots &  &  &  &  & \\
&  &  & \cos\alpha_{k} &  &  &  & \\
&  &  &  & \cos\alpha_{k} &  &  & \\
&  &  &  &  & 1 &  & \\
& \mathbf{0} &  &  &  &  & \ddots & \\
&  &  &  &  &  &  & 1
\end{pmatrix}
(Qx)\\
&  \geq\min_{j=1,...,k}\{\cos\alpha_{j},1\}(Qx)^{T}(Qx)=\min_{j=1,...,k}%
\{\cos\alpha_{j},1\}\left\vert Ox\right\vert ^{2}\\
&  =\min_{j=1,...,k}\{\cos\alpha_{j},1\}\left\vert x\right\vert ^{2},
\end{align*}
and subsequently%
\begin{align*}
&  \frac{d}{dt}\int_{%
\mathbb{R}
^{n}}u\left(  x\cdot P^{-1}x\right)  dx\\
&  \leq2Tr(P^{-1})M-\frac{\chi\min_{j=1,...,k}\{\cos\alpha_{j},1\}}{n\left\vert
B_{1}(0)\right\vert }\int_{%
\mathbb{R}
^{n}\times%
\mathbb{R}
^{n}}\frac{1}{\left\vert x-y\right\vert ^{n-2}}u(x,t)u(y,t)dydx.
\end{align*}
To simplify the last inequality, we invoke a result from \cite[Lemma
3.2.]{Biler}, which states that for any nonnegative function $f\in L^{1}(%
\mathbb{R}
^{n},(1+\left\vert x\right\vert ^{2})dx)$, the moment $m=\int_{%
\mathbb{R}
^{n}}f(x)\left\vert x\right\vert ^{2}dx$, the mass $M=\int_{%
\mathbb{R}
^{n}}f(x)dx$ and the integral $J:=\int_{%
\mathbb{R}
^{n}\times%
\mathbb{R}
^{n}}f(x)f(y)\left\vert x-y\right\vert ^{2-n}dydx,$ satisfy the inequality
$M^{\frac{n}{2}+1}\leq J(2m)^{\frac{n}{2}-1}$.

Therefore%
\begin{equation}
\int_{%
\mathbb{R}
^{n}\times%
\mathbb{R}
^{n}}\frac{1}{\left\vert x-y\right\vert ^{n-2}}u(x,t)u(y,t)dydx\geq
M^{\frac{n}{2}+1}\left(  2\int_{%
\mathbb{R}
^{n}}u\left\vert x\right\vert ^{2}dx\right)  ^{1-\frac{n}{2}}. \label{E1}%
\end{equation}
and the functions $w(t):=\int u(x,t)\left(  x\cdot P^{-1}x\right)  dx$ and
$m(t):=\int u(x,t)\left\vert x\right\vert ^{2}dx$ satisfy%
\begin{equation}
\frac{d}{dt}w(t)\leq2Tr(P^{-1})M-\frac{2^{1-\frac{n}{2}}\chi\min
_{j=1,...,k}\{\cos\alpha_{j},1\}M^{\frac{n}{2}+1}}{n\left\vert B_{1}%
(0)\right\vert }\left(  m(t)\right)  ^{1-\frac{n}{2}}. \label{E}%
\end{equation}
Let us denote the minimum and maximum eigenvalues of $P^{-1}$ by
$\lambda_{\min}$ and $\lambda_{\max},$ respectively. A standard result (cf.
\cite[Theorem 4.2.2.]{Roger}) asserts $\lambda_{\min}\left\vert x\right\vert
^{2}\leq x^{T}P^{-1}x\leq\lambda_{\max}\left\vert x\right\vert ^{2}$ for all
$x\in%
\mathbb{R}
^{n},$ yielding $\lambda_{\min}m(t)\leq w(t)\leq\lambda_{\max}m(t)$ for all
$t\geq0,$ and%
\begin{equation}
\frac{d}{dt}w(t)\leq2Tr(P^{-1})M-\frac{2^{1-\frac{n}{2}}\chi\min
_{j=1,...,k}\{\cos\alpha_{j},1\}M^{\frac{n}{2}+1}}{n\left\vert B_{1}%
(0)\right\vert }\left(  \lambda_{\min}\right)  ^{\frac{n}{2}-1}\left(
w(t)\right)  ^{1-\frac{n}{2}}. \label{E2}%
\end{equation}
This reads as the differential inequality%
\begin{align}
\frac{2}{n}\frac{d}{dt}w^{n/2}  &  \leq2Tr(P^{-1})Mw^{\frac{n}{2}-1}%
-\frac{2^{1-\frac{n}{2}}\chi\min_{j=1,...,k}\{\cos\alpha_{j},1\}M^{\frac{n}%
{2}+1}}{n\left\vert B_{1}(0)\right\vert }\left(  \lambda_{\min}\right)
^{\frac{n}{2}-1}\nonumber\\
&  =:f(w). \label{blow}%
\end{align}
We now introduce the condition on the initial data $f(w(0))<0.$ Since $f$ \ is
an increasing function of $w,$ the condition $f(w(0))<0$ implies that the
right-hand side of (\ref{blow}) is always negative and bounded away from zero.
We conclude that the right hand side is always negative and bounded away from
zero$.$ This leads to the conclusion that the function $w$ decreases and
assumes negative values in a finite time, contradicting the existence of a
global in time nonnegative solution. Finally, observing the inequality
$w(t)\leq\lambda_{\max}m(t),$ we obtain
\begin{align*}
f(w)  &  \leq2Tr(P^{-1})M\lambda_{\max}^{\frac{n}{2}-1}m^{\frac{n}{2}-1}%
-\frac{2^{1-\frac{n}{2}}\chi\min_{j=1,...,k}\{\cos\alpha_{j},1\}M^{\frac{n}%
{2}+1}}{n\left\vert B_{1}(0)\right\vert }\left(  \lambda_{\min}\right)
^{\frac{n}{2}-1}\\
&  =:h(m).
\end{align*}
Hence the condition on the initial moment $h(m(0)))<0$ or equivalently
\[
\int_{%
\mathbb{R}
^{n}}u(x,0)\left\vert x\right\vert ^{2}dx\leq\left(  \frac{2^{1-\frac{n}{2}%
}\chi\min_{j=1,...,k}\{\cos\alpha_{j},1\}}{2Tr(P^{-1}%
)\lambda_{\max}^{\frac{n}{2}-1}n\left\vert B_{1}(0)\right\vert }\left(
\lambda_{\min}\right)  ^{\frac{n}{2}-1}\right)  ^{\frac{2}{n-2}}M^{\frac{n}{n-2}},%
\] implies that $T_{\max}<\infty.$
\end{proof}

\begin{remark}
    For all \( M > 0 \) (even arbitrarily small), there exists an initial data \( u_0 \) with mass \( M \) such that the condition (\ref{SmallMoment}) is satisfied. Indeed, it is sufficient to consider non-negative, smooth, compactly supported data \( u_0 \) with mass \( M \) and second moment \( m_0 \). By rescaling it with \( \varepsilon^{-n} u_0\left(\frac{x}{\varepsilon}\right) \) for a sufficiently small \( \varepsilon > 0 \) (specifically, \( \varepsilon^2 \leq \frac{C_{Bl} M^{\frac{n}{n-2}}}{m_0} \)), the desired condition is achieved. In other words, blow-up is still possible for arbitrarily small initial mass, which contrasts with the two-dimensional case \cite{Cuentas}. 
\end{remark}

\section{\textbf{Global existence}}

\begin{theorem}
[Global existence]\label{GlobalHigh}Let $A:=(a_{ij})_{i,j=1,\ldots,n}\in M_{n}(%
\mathbb{R}
)$ be a matrix with constant components$.$ Then, there exists $\delta>0$ with
the property that if $\left\Vert u_{0}\right\Vert _{L^{\frac{n}{2}}}\leq\delta,$ for any non-negative $u_{0}\in
BUC(\mathbb{R}^{n})\cap L^{1}(\mathbb{R}^{n})$, the solution $u$ of the system
(\ref{Problem Matrix A nxn}) is global and for some constant $C>0,$ we have
that $\left\Vert u(\cdot,t)\right\Vert _{L^{\infty}(\mathbb{R}^{n})}\leq C$
for all $t>0.$
\end{theorem}

\begin{proof}
By multiplying the equation for $u$ by $u^{p-1}$ and integrating over
$\mathbb{R}^{n}$, we derive%
\begin{align*}
&  \frac{1}{p}\frac{d}{dt}\int_{\mathbb{R}^{n}}\left\vert u(x,t)\right\vert
^{p}dx\\
&  =-\frac{4(p-1)}{p^{2}}\int_{\mathbb{R}^{n}}\left\vert \nabla\left(
u^{p/2}\right)  \right\vert ^{2}dx-\frac{\chi(p-1)}{p}\int_{\mathbb{R}^{n}%
}u^{p}\left(  \nabla\cdot A\nabla v\right)  dx.
\end{align*}
Applying H\"{o}lder's inequality, we find%
\begin{align*}
&  \left\Vert u^{p}\left(  \nabla\cdot A\nabla v\right)  \right\Vert
_{L^{1}(\mathbb{R}^{n})}\\
&  \leq\left\Vert u\right\Vert _{L^{p+1}(\mathbb{R}^{n})}^{p}\left\Vert
\nabla\cdot A\nabla v\right\Vert _{L^{p+1}(\mathbb{R}^{n})}\leq\left\Vert
u\right\Vert _{L^{p+1}(\mathbb{R}^{n})}^{p}\sum_{_{i,j=1,n}}\left\vert
a_{ij}\right\vert \left\Vert \partial_{ij}\mathbf{K_n}\ast u\right\Vert
_{L^{p+1}(\mathbb{R}^{n})}\\
&  \leq\left\Vert u\right\Vert _{L^{p+1}(\mathbb{R}^{n})}^{p}\left\Vert
A\right\Vert _{\max}\sum_{_{i,j=1,n}}\left\Vert \partial_{ij}\mathbf{K_n}\ast
u\right\Vert _{L^{p+1}(\mathbb{R}^{n}).}%
\end{align*}
Now, we recall the following Calder\'{o}n--Zygmund inequality (See for
instance \cite[Section 6.4.2.]{Giga}): For all $g\in L^{q}(\mathbb{R}^{n}),$
there exist a constant $C_{CZI}^{(q,n)}=C(q,n),1<q<\infty,$ such that%
\begin{equation}
\left\Vert \partial_{ij}\mathbf{K_n}\ast g\right\Vert _{L^{q}(\mathbb{R}^{n}%
)}\leq C_{CZI}^{(q,n)}\left\Vert g\right\Vert _{L^{q}(\mathbb{R}^{n})},\text{
}i,j=1,2, \label{Calderon-Zygmund inequality}%
\end{equation}
Taking $g=u$ and $q=p+1,$ we deduce%
\[
\left\Vert u^{p}\left(  \nabla\cdot A\nabla v\right)  \right\Vert
_{L^{1}(\mathbb{R}^{n})}\leq4\left\Vert A\right\Vert _{\max}C_{CZI}%
^{(p+1,n)}\left\Vert u\right\Vert _{L^{p+1}(\mathbb{R}^{n})}^{p+1}.
\]
This leads to%
\begin{align}
&  \frac{1}{\left(  p-1\right)  }\frac{d}{dt}\int_{\mathbb{R}^{n}}\left\vert
u(x,t)\right\vert ^{p}dx\nonumber\\
&  \leq-\frac{4}{p}\int_{\mathbb{R}^{n}}\left\vert \nabla\left(
u^{p/2}\right)  \right\vert ^{2}dx+4\chi\left\Vert A\right\Vert _{\max}%
C_{CZI}^{(p+1,n)}\int_{\mathbb{R}^{n}}u^{p+1}dx. \label{f1}%
\end{align}
Applying the Gagliardo-Nirenberg-Sobolev inequality, we obtain that for any
$\frac{n}{2}\leq p+1\leq\frac{pn}{n-2}$%
\begin{align}
\int_{\mathbb{R}^{n}}u^{p+1}dx  &  \leq\left\Vert u\right\Vert _{L^{\frac
{n}{2}}}\left\Vert u\right\Vert _{\frac{pn}{n-2}}^{p}=\left\Vert u\right\Vert
_{L^{\frac{n}{2}}}\left\Vert u^{\frac{p}{2}}\right\Vert _{L^{\frac{2n}{n-2}}%
}^{2}\nonumber\\
&  \leq C_{GNS}^{2}\left\Vert u\right\Vert _{L^{\frac{n}{2}}}\int
_{\mathbb{R}^{n}}\left\vert \nabla\left(  u^{p/2}\right)  \right\vert ^{2}dx.
\label{f2}%
\end{align}
Combining (\ref{f1}) and (\ref{f2}), we get for any $p\geq\max\{1,\frac{n}%
{2}-1\}$%
\begin{align}
&  \frac{1}{\left(  p-1\right)  }\frac{d}{dt}\int_{\mathbb{R}^{n}}\left\vert
u(x,t)\right\vert ^{p}dx\nonumber\\
&  \leq\left(  4\chi\left\Vert A\right\Vert _{\max}C_{CZI}^{(p+1,n)}%
C_{GNS}^{2}\left\Vert u\right\Vert _{L^{\frac{n}{2}}}-\frac{4}{p}\right)
\int_{%
\mathbb{R}
^{n}}\left\vert \nabla\sqrt{u}\right\vert ^{2}dx. \label{f3}%
\end{align}
Notice that for $p=\frac{n}{2}$ in (\ref{f3})$,$ the inequality $4\chi
\left\Vert A\right\Vert _{\max}C_{CZI}^{(n/2+1,n)}C_{GNS}^{2}\left\Vert
u_{0}\right\Vert _{L^{\frac{n}{2}}}-\frac{8}{n}\leq0$ implies that $\left\Vert
u\right\Vert _{L^{\frac{n}{2}}}$ decreases for $t\in\left(  0,T_{\max}\right)
.$ As a consequence the condition%
\[
\left\Vert u_{0}\right\Vert _{L^{\frac{n}{2}}}\leq\frac{1}{\chi\left\Vert
A\right\Vert _{\max}C_{GNS}^{2}}\min\left\{  \frac{2}{nC_{CZI}^{(n/2+1,n)}%
},\frac{1}{pC_{CZI}^{(p+1,n)}}\right\}  =:\delta(p,n),
\]
for $p\geq\max\{1,\frac{n}{2}-1\}$ implies that the function $\int
_{\mathbb{R}^{n}}\left\vert u(x,t)\right\vert ^{p}dx$ decreases for
$t\in\left(  0,T_{\max}\right).$ \newline We fix any $q>n,$ and let
$\delta:=\delta(q,n).$ Then, assuming that $\left\Vert u_{0}\right\Vert
_{L^{\frac{n}{2}}}\leq\delta,$ we obtain from (\ref{f3}) that there exists $c_{1}>0$ such that
\begin{equation}
\left\Vert u(\cdot,t)\right\Vert _{L^{q}(\mathbb{R}^{n})}\leq c_{1}\text{ for
all }t\in(0,T_{\max}).\label{PB}%
\end{equation}
We recall now the following $L^{q}-L^{p}$ estimates of heat semigroup
$e^{t\Delta}$. For any $1\leq q\leq p\leq\infty,$ there holds%
\begin{align}
\left\Vert e^{t\Delta}f\right\Vert _{L^{p}(\mathbb{R}^{n})} &  \leq(4\pi
t)^{\frac{n}{2}(\frac{1}{p}-\frac{1}{q})}\left\Vert f\right\Vert
_{L^{q}(\mathbb{R}^{n})},\label{HS1}\\
\left\Vert \nabla\cdot e^{t\Delta}F\right\Vert _{L^{p}(\mathbb{R}^{n})} &
\leq Ct^{-\frac{1}{2}+\frac{n}{2}\left(  \frac{1}{p}-\frac{1}{q}\right)
}\left\Vert F\right\Vert _{L^{q}(\mathbb{R}^{n})},\label{HS2}%
\end{align}
where $C=C(p,q,n)$ is a constant depending only on $p,q$ and $n.$ These
inequalities are a consequences of Young's inequality for the convolution (For
example, see \cite[Subsection 4.1.2. p. 145]{Giga}).

Let us define%
\[
N(T):=\sup_{t\in(0,T)}\left\Vert u(\cdot,t)\right\Vert _{L^{\infty}%
(\mathbb{R}^{n})},\text{ for }T\in(0,T_{\max}).
\]
Using the Duhamel integral equation, we get%
\[
u(t)=e^{(t-t_{0})\Delta}u(t_{0})-\chi\int_{t_{0}}^{t}\nabla\cdot e^{\left(
t-s\right)  \Delta}\left(  u(s)A\nabla v(s)\right)  ds,
\]
with%
\[
t_{0}=\left\{
\begin{array}
[c]{cc}%
t-1, & \text{if }t\geq1,\\
0, & \text{if }0\leq t\leq1.
\end{array}
\right.
\]
By (\ref{HS1}) and (\ref{HS2}), we have that%
\begin{align*}
& \left\Vert u(x,t)\right\Vert _{L^{\infty}(\mathbb{R}^{n})}\\
& \leq\left\Vert e^{(t-t_{0})\Delta}u(t_{0})\right\Vert _{L^{\infty
}(\mathbb{R}^{n})}+\chi\left\Vert \int_{t_{0}}^{t}\nabla\cdot e^{\left(
t-s\right)  \Delta}\left(  u(s)A\nabla v(s)\right)  ds\right\Vert _{L^{\infty
}(\mathbb{R}^{n})}\\
& \leq(4\pi(t-t_{0}))^{\frac{-n}{2q}}\left\Vert u(t_{0})\right\Vert
_{L^{q}(\mathbb{R}^{n})}\\
& +c_{2}\int_{t_{0}}^{t}(t-s)^{-\frac{1}{2}-\frac{n}{2q}}\left\Vert \left\vert
u(s)A\nabla v(s)\right\vert \right\Vert _{L^{q}(\mathbb{R}^{n})}ds\\
& \leq(4\pi(t-t_{0}))^{\frac{-n}{2q}}\left\Vert u(t_{0})\right\Vert
_{L^{q}(\mathbb{R}^{n})}\\
& +c_{2}\left\Vert A\right\Vert _{\max}\int_{t_{0}}^{t}(t-s)^{-\frac{1}%
{2}-\frac{n}{2q}}\left\Vert u(s)\right\Vert _{L^{q}(\mathbb{R}^{n})}\left\Vert
\nabla v(s)\right\Vert _{L^{\infty}(\mathbb{R}^{n})}ds,
\end{align*}
for all $t\in\left(  0,T\right)  .$ Notice that for any $\gamma>0,$ we have
that%
\begin{align}
& \left\vert \nabla v(x,s)\right\vert \nonumber\\
& =\left\vert \nabla\mathbf{K_n}(x)\ast u(x,s)\right\vert \nonumber\\
& =\left\vert \frac{-1}{n\left\vert B_{1}(0)\right\vert }\int_{%
\mathbb{R}
^{n}}\frac{x-y}{\left\vert x-y\right\vert ^{n}}u(y,s)dy\right\vert
\nonumber\\
& \leq\frac{1}{n\left\vert B_{1}(0)\right\vert }\left(  \int_{\left\vert
x-y\right\vert \leq\gamma}\frac{u(y,s)}{\left\vert x-y\right\vert ^{n-1}%
}dy+\int_{\left\vert x-y\right\vert >\gamma}\frac{u(y,s)}{\left\vert
x-y\right\vert ^{n-1}}dy\right)  \nonumber\\
& \leq\frac{\left\Vert u(x,s)\right\Vert _{L^{\infty}(\mathbb{R}^{n})}%
}{n\left\vert B_{1}(0)\right\vert }\int_{\left\vert z\right\vert \leq\gamma
}\left\vert z\right\vert ^{1-n}dz+\frac{\gamma^{1-n}\left\Vert
u(x,s)\right\Vert _{L^{1}(\mathbb{R}^{n})}}{n\left\vert B_{1}(0)\right\vert
}\nonumber\\
& =\gamma\left\Vert u(x,s)\right\Vert _{L^{\infty}(\mathbb{R}^{n})}%
+\frac{\gamma^{1-n}M}{n\left\vert B_{1}(0)\right\vert }.\label{VB}%
\end{align}
Therefore, from (\ref{PB}) and (\ref{VB})%
\begin{align}
& \left\Vert u(x,t)\right\Vert _{L^{\infty}(\mathbb{R}^{n})}\nonumber\\
& \leq c_{1}(4\pi)^{\frac{-n}{2q}}\nonumber\\
& +c_{1}c_{2}\left\Vert A\right\Vert _{\max}\left(  \gamma N(T)+\frac
{\gamma^{1-n}M}{n\left\vert B_{1}(0)\right\vert }\right)  \int_{t_{0}}%
^{t}(t-s)^{-\frac{1}{2}-\frac{n}{2q}}ds,\label{IN}%
\end{align}
for all $t\in\left(  0,T\right)  .$ Note that%
\begin{align}
\int_{t_{0}}^{t}(t-s)^{-\frac{1}{2}-\frac{n}{2q}}ds  & =\int_{0}^{t-t_{0}}%
\tau^{-\frac{1}{2}-\frac{n}{2q}}d\tau\nonumber\\
& \leq\int_{0}^{1}\tau^{-\frac{1}{2}-\frac{n}{2q}}d\tau\nonumber\\
& =\frac{2q}{q-n}.\label{IN2}%
\end{align}
Taking%
\[
\frac{2qc_{1}c_{2}\left\Vert A\right\Vert _{\max}\gamma}{q-n}=\frac{1}{2},
\]
we conclude from (\ref{IN}) and (\ref{IN2}) that there exists a constant
\[
c_{3}:=c_{1}(4\pi)^{\frac{-n}{2q}}+c_{1}c_{2}\left\Vert A\right\Vert _{\max
}\frac{\gamma^{1-n}M}{n\left\vert B_{1}(0)\right\vert }>0,
\]
such that
\[
N(T)\leq\frac{1}{2}N(T)+c_{3}\text{ for all }T<T_{\max},
\]
and hence%
\[
\left\Vert u(\cdot,t)\right\Vert _{L^{\infty}(\mathbb{R}^{n})}\leq2c_{3}\text{
for all }t\in\left(  0,T_{\max}\right)  ,
\]
as $T\in\left(  0,T_{\max}\right)  $ was arbitrary. Taking into account the
extensibility criterion in Preposition \ref{LocalExistence nxn}, the last
inequality implies global existence. \bigskip
\end{proof}

\begin{remark}
Applying the inequality that compares the \(L^{\frac{n}{2}}\)-norm, the mass \(M\), and the second moment \(m_0\) of a non-negative function \(u_0\) (See \cite[Remark 2.6]{BilerKarch}):
\begin{equation}
\|u_0\|_{L^{\frac{n}{2}}} \geq C_n M \left( \frac{M}{m_0} \right)^{\frac{n-2}{2}}, \tag{6.21}
\end{equation}
where \(C_n = C(n)\) is a constant depending only on \(n\), we find that the condition (\ref{SmallMoment}) in Theorem \ref{blow-up nxn} implies: 
\begin{equation*}
    \|u_0\|_{L^{\frac{n}{2}}} \geq C_n (C_{\text{Bl}})^{\frac{2-n}{2}}.
\end{equation*}
Conversely, the smallness assumption on \(\|u_0\|_{L^{\frac{n}{2}}}\) in Theorem \ref{GlobalHigh} implies: 
\begin{equation*}
    m_0\geq \left(\frac{C_n}{\delta}\right)^{\frac{2}{n-2}}M^{\frac{n}{n-2}},
\end{equation*}
which shows the compatibility of both results.
\end{remark}

\end{document}